\pdfoutput=1
\documentclass[11pt,a4paper]{article}
\usepackage[USenglish]{babel}
\usepackage[T1]{fontenc}
\usepackage[utf8]{inputenc}
\usepackage{titlesec}
\usepackage{csquotes}
\usepackage{enumitem}
\usepackage{appendix}
\usepackage[giveninits,urldate=comp,sortcites,maxbibnames=10,url=false]{biblatex}
\usepackage[colorlinks=true,linkcolor=blue,bookmarksnumbered=true,bookmarksdepth=3,pdftitle=Representation\ Formula\ for\ Viscosity\ Solutions\ to\ a\ class\ of\ Nonlinear\ Parabolic\ PDEs,pdfauthor=Marco\ Pozza]{hyperref}

\usepackage{amssymb}
\usepackage{amsmath}
\usepackage{mathtools}
\usepackage{amsthm}
\usepackage{thmtools}
\usepackage{thm-restate}
\usepackage{dsfont}

\usepackage[noabbrev]{cleveref}

\titleformat*{\subsubsection}{\itshape}

\theoremstyle{plain}
\newtheorem{theo}{Theorem}[section]
\crefname{theo}{theorem}{theorems}
\crefname{dpp}{dynamic programming principle}{}
\newtheorem{prop}[theo]{Proposition}

\newtheorem{lem}[theo]{Lemma}
\theoremstyle{definition}
\newtheorem{defin}[theo]{Definition}

\newtheorem{prob}[theo]{Problem}
\newtheorem{assum}[theo]{Assumptions}
\crefname{assum}{assumptions}{assumptions}
\theoremstyle{remark}
\newtheorem{rem}[theo]{Remark}
\crefname{rem}{remark}{remarks}

\AtBeginEnvironment{appendices}{\crefalias{section}{appendix}}

\newlist{myenum}{enumerate}{1}
\setlist[myenum]{label={\bf\roman*)}}
\crefname{myenumi}{item}{items}

\DeclareMathOperator{\tr}{tr}
\DeclareMathOperator*{\esss}{ess\,sup}
\DeclareMathOperator*{\essi}{ess\,inf}

\newcommand{\Acal}{\mathcal A}
\newcommand{\Bcal}{\mathcal B}
\newcommand{\Eds}{\mathds E}
\newcommand{\Fcal}{\mathcal F}
\newcommand{\Ical}{\mathcal I}
\newcommand{\Jcal}{\mathcal J}
\newcommand{\Nds}{\mathds N}
\newcommand{\Pcal}{\mathcal P}
\newcommand{\Pds}{\mathds P}
\newcommand{\Rds}{\mathds R}
\newcommand{\Sds}{\mathds S}

\title{Representation Formula for Viscosity Solutions to a class of Nonlinear Parabolic PDEs}
\author{Marco Pozza\thanks{Dipartimento di Matematica ``G. Castelnuovo'', Sapienza Università di Roma. Piazzale Aldo Moro 5, 00185 Roma, Italy. {\tt marco.pozza@uniroma1.it}}}
\date{}

\begin{document}

\maketitle

\begin{abstract}
We provide a representation formula for viscosity solutions to a class of nonlinear second order parabolic PDEs given as a sup--envelope function. This is done through a dynamic programming principle derived from \cite{art:pathpeng}. The formula can be seen as a nonlinear extension of the Feynman--Kac formula and is based on the backward stochastic differential equations theory.
\end{abstract}


\medskip


\section{Introduction}

It is well known that viscosity solutions were conceived by Crandall and Lions (1982) in the framework of optimal control theory. The goal was to show well posedness of Hamilton--Jacobi--Bellman equations in the whole space, and to prove, via dynamic programming principle, the value function of a suitable optimal control problem being the unique solution.

When trying to extend viscosity methods to the analysis of second order \emph{parabolic partial differential equations}, PDEs for short, and getting representation formulas, it appeared clear that some stochastic dynamics must be brought into play. Not surprisingly, this has been first done for stochastic control models. The hard work of a generation of mathematicians, \cite{hjlions2,hjlions1,bensoussa82,krylovbook,flemingrishel75,nisio15} among the others, allowed making effective dynamic programming approach to stochastic control problems.

Prompted by this body of investigations, a stream of research arose in the probabilistic community ultimately leading to the theory of \emph{backward stochastic differential equations}, BSDEs for short, which was introduced by Pardoux and Peng in \cite{Pardoux_Peng_1990} (1990). Since then, it has attracted a great interest due to its connections with mathematical finance and PDEs, as well as with stochastic control. This theory has been in particular used to extend the classical Feynman--Kac formula, which establishes a link between linear parabolic PDEs and \emph{stochastic differential equations}, SDEs for short, to semilinear and quasilinear equations, see for example \cite{Delarue_2002,Delarue_Guatteri_2006,Ma_Yong_2007,pardouxpeng92}. See also \cite{pardouxbook} for a rather complete overview of the semilinear case.

For sake of clarity, let us consider the following semilinear parabolic PDE problem coupled with final conditions,
\begin{equation}\label{eq:exparsemlin}
\begin{cases}
\begin{aligned}
\frac12\left\langle\sigma\sigma^\dag,D^2_xu\right\rangle(t,x)+(\nabla_xub)(t,x)&\\
+\partial_tu(t,x)+f(t,x,u,\nabla_xu\sigma)&=0
\end{aligned}
&t\in[0,T],x\in\Rds^N,\\
u(T,x)=g(x),&x\in\Rds^N,
\end{cases}
\end{equation}
then its viscosity solution can be written as $u(t,x)=\Eds\left(Y^{t,x}_t\right)$, where $Y$ is given by the following system, called \emph{forward backward stochastic differential equation} or FBSDE in short, which is made in turn of two equations, the first one is a SDE, and the second one a \emph{backward stochastic differential equation}, BSDE for short, depending on the first one:
\[
\left\{\begin{aligned}
\!X^{t,x}_s\!=&x+\int_t^s\sigma\left(r,X^{t,x}_r\right)dW_r+\int_t^sb\left(r,X^{t,x}_r\right)dr,\\
\!Y^{t,x}_s\!=&g\left(X^{t,x}_T\right)\!+\!\!\int_s^T\!\!f\!\left(r,X_r^{t,x},Y_r^{t,x},Z_r^{t,x}\right)\!dr\!-\!\!\int_s^T\!\!Z_r^{t,x}dW_r,
\end{aligned}\right.s\!\in\![t,T],x\!\in\!\Rds^N.
\]
As can be intuitively seen, the SDE takes care of the linear operator defined by $\sigma$ and $b$, also called the \emph{infinitesimal generators} of the SDE, while the BSDE depends on $f$ and $g$. In other words, this extension of the Feynman--Kac formula basically does not modify the treatment of the second order linear operator with respect to the completely linear case.

Subsequently, Peng introduced in \cite{math/0601035} (2006) the notion of \emph{$G$--expectation}, a nonlinear expectation generated by a fully nonlinear second order operator $G$ via its viscosity solutions. This work has originated an active research field, with relevant applications to Mathematical Finance.

Peng has improved this theory in several papers and has given a comprehensive account of it in the book \cite{book:peng}, where he highlights the role of the so--called \emph{sublinear expectations}, namely $G$--expectations generated by sublinear operators given as a sup--envelope function. Finally in \cite{art:pathpeng}, Peng provides representation formulas for viscosity solutions using these expectations. More precisely, given a sublinear operators $G$ and the $G$--heat equation
\begin{equation}\label{eq:gheat}
\begin{cases}
\partial_tu(t,x)+G\left(D_x^2u(t,x)\right)=0,&t\in[0,T],x\in\Rds^N,\\
u(T,x)=g(x),&x\in\Rds^N.
\end{cases}
\end{equation}
he represents the viscosity solution as
\[
u(t,x)=\sup_{\sigma\in\Acal}\Eds\left(g\left(x+\int_t^T\sigma_sdW_s\right)\right),
\]
where $\Acal$ is a family of stochastic process associated to $G$ and $W$ is a Brownian motion.
The key to prove it is a \emph{dynamic programming principle} that we will illustrate in the paper, see \cref{dynprincsec}. We point out that this formula has been established in \cite{hjlions2} using a different method and that the novelty with respect to the Feynman--Kac formula is essentially given by the sublinearity of the operator $G$.

The purpose of this article is to apply a generalized version of the dynamic programming principle of \cite{art:pathpeng} in order to give representation formulas of solutions to PDE problems of the type
\begin{equation}\label{eq:parsubopgen}
\begin{cases}
\partial_tu(t,x)+F\left(t,x,u,\nabla_xu,D^2_xu\right)=0,&t\in[0,T],x\in\Rds^N,\\
u(T,x)=g(x),&x\in\Rds^N,
\end{cases}
\end{equation}
where $F$ is a sup--envelope function of semilinear operators as in \eqref{eq:exparsemlin}. This problem is clearly a blend between \eqref{eq:exparsemlin} and \eqref{eq:gheat} which generalize both.

This is hopefully just a first step to further extend the Feynman--Kac formula to fully nonlinear problems. In this regard the method here employed has been successfully used in \cite{pozza21} to provide a representation formula for a particular Dirichlet type problem. We also point out that there is close connection between this method and second order BSDEs, 2BSDEs for short. 2BSDEs were introduced by Cheridito, Soner, Touzi and Victoir in \cite{int2bsde} (2007). Then, in 2011, Soner, Touzi and Zhang \cite{sonertouzizhang11} provided a complete theory of existence and uniqueness for 2BSDEs under Lipschitz conditions. In those papers is also analyzed the connection between 2BSDEs and fully nonlinear PDEs. Among the subsequent developments of this theory we cite \cite{possamaitanzhou18,sonertouzizhang11,ekrentouzizhang16,ekrentouzizhang16_2,matoussipossamaizhou13,kazi-tanipossamaizhou15} and in particular \cite{possamai13}, which performs its analysis replacing the Lipschitz condition on $y$ with monotonicity as we do here.

This paper is organized as follows: in \cref{setprob} we introduce the problem along with some preliminary facts and in \cref{dynprincsec} we develop a dynamic programming principle which is the core of our theory. Then, in \cref{repform}, we perform the essential part of our analysis, and obtain in this way our main results. In \cref{2bsdesec} we summarily analyze the connection between 2BSDEs and our representation formula, while \cref{sdeprel} briefly gives some probability results about stochastic differential equations needed to our analysis. Finally in \cref{compsec} is proved a comparison result for parabolic problems.

We proceed setting the notation used in the paper.

\subsubsection*{Notation}
\pdfbookmark[3]{Notation}{Notation}

We will work on the filtered probability space $\left(\Omega,\Fcal,\{\Fcal_t\}_{t\in[0,\infty)},\Pds\right)$,
\begin{itemize}
\item $\Fcal$ is a complete $\sigma$--algebra on $\Omega$;
\item the stochastic process $\{W_t\}_{t\in[0,\infty)}$ will denote the $N$ dimensional Brownian motion under $\Pds$;
\item $\{\Fcal_t\}_{t\in[0,\infty)}$ is the filtration defined by $\{W_t\}_{t\in[0,\infty)}$ which respects the usual condition of completeness and right continuity;
\item $\{W_s^t\}_{s\in[t,\infty)}:=\{W_s-W_t\}_{s\in[t,\infty)}$ is a Brownian motion independent from $\{W_s\}_{s\in[0,t]}$ by the \emph{strong Markov property};
\item $\{\Fcal^t_s\}_{s\in[t,\infty)}$ is the filtration generated by $\{W_s^t\}_{s\in[t,\infty)}$ which we assume respect the usual condition and is independent from $\Fcal_t$;
\item we will say that a stochastic process $\{H_t\}_{t\in[0,\infty)}$ is adapted if $H_t$ is $\Fcal_t$--measurable for any $t\in[0,\infty)$;
\item we will say that a stochastic process $\{H_t\}_{t\in[0,\infty)}$ is progressively measurable, or simply progressive, if, for any $T\in[0,\infty)$, the application that to any $(t,\omega)\in[0,T]\times\Omega$ associate $H_t(\omega)$ is measurable for the $\sigma$--algebra $\Bcal([0,T])\times\Fcal_T$;
\item a function on $\Rds$ is called \emph{cadlag} if is right continuous and has left limit everywhere;
\item a cadlag (in time) process is progressive if and only if is adapted;
\item $B_\delta(x)$ will denote an open ball centered in $x$ with radius $\delta$;
\item if $A\in\Rds^{N\times M}$ then $A^\dag$ will denote its transpose;
\item (Frobenius product) if $A,B\in\Rds^{N\times M}$ then $\langle A,B\rangle:=\tr\left(AB^\dag\right)$ and $|A|$ is the norm $\sqrt{\langle A,A\rangle}=\sqrt{\sum\limits_{i=1}^N\sum\limits_{j=1}^MA_{i,j}^2}$;
\item $\Sds^N$ is the space of all $N\times N$ real valued symmetric matrices and $\Sds^N_+$ is the subset of $\Sds^N$ made up by the definite positive matrices.
\end{itemize}

\noindent\textbf{Acknowledgement:} We would like to thank an anonymous referee to a previous version for suggesting us the comparison \cref{comptheo} as well as other comments which helped us improve this paper.

\section{Setting of the Problem}\label{setprob}

Here we introduce the problem, which is object of investigation, as well as basic preliminary facts.

We start giving the assumption of the PDE operator which we will analyze in this paper.

\begin{assum}\label{probass}
We fix a $T>0$, hereafter called \emph{terminal time}, the constants $\ell\ge0$, $\lambda\ge0$, $\mu\in\Rds$ and let $K_F$ be a set made up by functions
\[
(b,\sigma):[0,T]\times\Rds^N\longrightarrow\Rds^N\times\Rds^{N\times M}
\]
and which satisfies, for any $(t,x)\in[0,T]\times\Rds^N$, the following conditions:
\begin{myenum}[series=probass]
\item\label{paloceb} (local equiboundedness) $\sup\limits_{(b,\sigma)\in K_F}|b(t,x)|+|\sigma(t,x)|<\infty$;
\item\label{paecont} (equicontinuity) for each $\varepsilon>0$ there is a neighborhood $U$ of $(t,x)$ such that $|b(t,x)-b(t',x')|+|\sigma(t,x)-\sigma(t',x')|<\varepsilon$ for any $(t',x')\in U$ and $(b,\sigma)\in K_F$;
\item\label{paelip} (equiLipschitz continuity) for any $(b,\sigma)\in K_F$ and $t\in[0,T]$
\[
|b(t,x)-b(t,x')|\le\ell|x-x'|,\qquad|\sigma(t,x)-\sigma(t,x')|\le\ell|x-x'|.
\]
\end{myenum}
Now let
\[
f:[0,T]\times\Rds^N\times\Rds^N\times\Rds^{N\times M}\times\Rds\times\Rds^N\longrightarrow\Rds
\]
be a continuous function such that, for any $t\in[0,T]$, $x,x',q,q',z,z'\in\Rds^N$, $A,A'\in\Rds^{N\times M}$ and $y,y'\in\Rds$,
\begin{myenum}[resume=probass]
\item $|f(t,x,q,A,y,z)-f(t,x',q',A',y,z')|\le\ell(|x-x'|+|q-q'|+|A-A'|+|z-z'|)$;
\item $|f(t,x,q,A,y,z)|\le\ell(1+|x|+|q|+|A|+|y|+|z|)$;
\item $(y-y')(f(t,x,q,A,y,z)-f(t,x,q,A,y',z))\le\mu|y-y'|^2$;
\item fixed a $(t,y)\in[0,T]\times\Rds$, a compact set $K\subset\Rds^N\times\Rds^{N\times M}\times\Rds^N$ and a positive $\varepsilon$ there is a $\delta>0$ such that
\[
|f(t,x,q,A,y,z)-f(t',x,q,A,y',z)|<\varepsilon
\]
for any $(t',y')\in B_\delta((t,y))\cap[0,T]\times\Rds$ and $(q,A,z)\in K$.
\end{myenum}
In what follows we will focus on nonlinear operator of this type
\[
F(t,x,y,p,S):=\sup_{(b,\sigma)\in K_F}\frac12\left\langle\sigma\sigma^\dag(t,x),S\right\rangle+p^\dag b(t,x)+f_{(b,\sigma)}\left(t,x,y,p^\dag\sigma\right),
\]
where $f_{(b,\sigma)}\left(t,x,y,p^\dag\sigma\right):=f\left(t,x,b(t,x),\sigma(t,x),y,p^\dag\sigma(t,x)\right)$.
\end{assum}

It is a simple verification that $F$ is continuous. Furthermore we have that $F$ is an elliptic and convex (for the second order term) operator in the sense that, for any $t\in[0,T]$, $x\in\Rds^N$, $y\in\Rds$, $p\in\Rds^N$ and $S\in\Sds^N$,
\begin{myenum}
\item\emph{(Convexity)} for any $\delta\in[0,1]$
\[
F(t,x,y,p,\delta S+(1-\delta)S')\le\delta F(t,x,y,p,S)+(1-\delta)F(t,x,y,p,S');
\]
\item\label{Fellip}\emph{(Ellipticity)} there is a $\lambda\ge0$ such that, for any $S'\ge0$,
\[
F(t,x,y,p,S+S')-F(t,x,y,p,S)\ge\lambda|S'|.
\]
\end{myenum}

The $\lambda$ in \cref{Fellip} is actually the lowest eigenvalue of the matrix $\sigma\sigma^\dag$.

\begin{rem}\label{probass2}
We point out that \cref{paecont,paloceb} in \cref{probass} and the Arzelà--Ascoli theorem imply that $K_F$ is relatively compact in the compact--open topology and it is easy to see that its closure $\overline{K_F}$ satisfies \cref{paecont,paelip,paloceb}. Moreover the map
\[
(b,\sigma)\mapsto\frac12\left\langle\sigma\sigma^\dag(t,x),S\right\rangle+p^\dag b(t,x)+f_{(b,\sigma)}\left(t,x,y,p^\dag\sigma\right)
\]
is continuous, thus we have that
\[
F(t,x,y,p,S)=\max_{(b,\sigma)\in\overline{K_F}}\frac12\left\langle\sigma\sigma^\dag(t,x),S\right\rangle+p^\dag b(t,x)+f_{(b,\sigma)}\left(t,x,y,p^\dag\sigma\right).
\]
For this reason we will always assume that $K_F$ is closed.
\end{rem}

\begin{prob}\label{parprob}
Let $T$ be a terminal time, $F$ a sublinear operator as in \cref{probass} and $g:\Rds^N\to\Rds$ a continuous function such that, for any $x,x'\in\Rds^N$,
\begin{myenum}
\item $|g(x)-g(x')|\le\ell|x-x'|$;
\item $|g(x)|\le\ell(1+|x|)$,
\end{myenum}
where $\ell$ is the same constant in \cref{probass}.\\
Find the solution $u$ to the parabolic PDE
\[
\begin{cases}
\partial_tu(t,x)+F\left(t,x,u,\nabla_xu,D^2_xu\right)=0,&t\in(0,T),x\in\Rds^N,\\
u(T,x)=g(x),&x\in\Rds^N.
\end{cases}
\]
\end{prob}

The solutions to \cref{parprob} are to be intended in the viscosity sense, hence let us define what we mean with viscosity solution. For a detailed overview of the viscosity solution theory we refer to \cite{userguide}.

\begin{defin}
Given an upper semicontinuous function $u$ we say that a function $\varphi$ is a \emph{supertangent} to $u$ at $(t,x)$ if $(t,x)$ is a local maximizer of $u-\varphi$.\\
Similarly we say that a function $\psi$ is a \emph{subtangent} to a lower semicontinuous function $v$ at $(t,x)$ if $(t,x)$ is a local minimizer of $v-\psi$.
\end{defin}

\begin{defin}
An upper semicontinuous function $u$ is called a \emph{viscosity subsolution} to \cref{parprob} if, for any suitable $(t,x)$ and $C^{1,2}$ supertangent $\varphi$ to $u$ at $(t,x)$,
\[
\partial_t\varphi(t,x)+F\left(t,x,u,\nabla_x\varphi,D^2\varphi\right)\ge0.
\]
Similarly a lower semicontinuous function $v$ is called a \emph{viscosity supersolution} to \cref{parprob} if, for any suitable $(t,x)$ and $C^{1,2}$ subtangent $\psi$ to $v$ at $(t,x)$,
\[
\partial_t\psi(t,x)+F\left(t,x,v,\nabla_x\psi,D^2\psi\right)\le0.
\]
Finally a continuous function $u$ is called a \emph{viscosity solution} to \cref{parprob} if it is both a super and a subsolution to \cref{parprob}.
\end{defin}

For viscosity sub and supersolutions to \cref{parprob} holds a comparison result which we will prove in \cref{compsec}.

\begin{restatable}{theo}{comppar}\label{compparprob}
Let $u$ and $v$ be respectively a sub and a supersolution to \cref{parprob}. If $u(T,x)\le v(T,x)$ for any $x\in\Rds^N$ and both $u$ and $v$ have polynomial growth, then $u\le v$ on $(0,T]\times\Rds^N$.
\end{restatable}

Our method to obtain representation formulas relies on a \emph{dynamic programming principle}, which is based on a construction on a broader set than $K_F$. This set, which we call $\Acal_F$, is made up of the controls
\[
(b,\sigma):[0,T]\times\Omega\times\Rds^N\longrightarrow\Rds^N\times\Rds^{N\times M},
\]
where $\Omega$ is a filtered probability space as in the notation section. These controls are cadlag, i.e. right continuous with finite left limits everywhere, on $[0,T]$ and such that, for any $t\in[0,T]$ and $\omega\in\Omega$, there is a neighborhood $(t_1,t_2)$ of $t$ and a pair $(b',\sigma')\in K_F$ such that $(b,\sigma)|_{(t_1,t_2)}(\omega)\equiv(b',\sigma)'|_{(t_1,t_2)}$. We also require that, for any $x\in\Rds^N$ and $(b,\sigma)\in\Acal_F$, the process $\{(b,\sigma)(t,x)\}_{t\in[0,T]}$ is progressively measurable. $\Acal_F$ is trivially non empty, since it contains $K_F$. From the definition we have that the elements of $\Acal_F$ are uniformly Lipschitz continuous in $\Rds^N$ with Lipschitz constant $\ell$ and the eigenvalues of $(\sigma\sigma^\dag)(t,\omega,x)$ belong to $[\lambda,\infty)^N$ for any $t\in[0,T]$, $\omega\in\Omega$, $x\in\Rds^N$. In what follows we will usually omit the dependence on $\Omega$.

For any stopping time $\tau$ with value in $[0,T]$, an useful subset of $\Acal_F$, which we will use later, is $\Acal^\tau_F$, which consists of the $(b,\sigma)$ belonging to $\Acal_F$ such that $\{(b,\sigma)(\tau+t,x)\}_{t\in[0,\infty)}$ is progressive with respect to the filtration $\left\{\Fcal_t^\tau\right\}_{t\in[0,\infty)}$. Clearly $\Acal^0_F=\Acal_F$ and
\[
F(t,x,y,p,S)=\max_{(b,\sigma)\in\Acal_F^t}\frac12\left\langle\sigma\sigma^\dag(t,x),S\right\rangle+p^\dag b(t,x)+f_{(b,\sigma)}\left(t,x,y,p^\dag\sigma\right).
\]

\begin{rem}\label{aflp}
By \cref{probass2}, for any bounded subset $A$ of $\Rds^N$ and $p>0$ we have that
\[
\sup_{(b,\sigma)\in\Acal_F}\Eds\left(\int_0^T\int_A\left(|b(t,x)|^p+|\sigma(t,x)|^p\right)dxds\right)<\infty.
\]
Thus $\Acal_F$ is contained in $L^p\left([0,T]\times\Omega\times A;\Rds^N\times\Rds^{N\times M}\right)$ for any bounded subset $A$ of $\Rds^N$ and $p>0$.
\end{rem}

\section{Dynamic Programming Principle}\label{dynprincsec}

The scope of this \namecref{dynprincsec} is to provide the necessary instruments to prove the \cref{cauchydpp}, which will be used to derive representation formulas for the viscosity solution to \cref{parprob}. The dynamic programming principle is, in this contest, an instrument that permit us to break a stochastic trajectory in two or more part. In particular for the problem \eqref{eq:exparsemlin} it means that
\begin{multline*}
\Eds\left(u\left(s,X^{t,x}_s\right)+\int_s^Tf\left(r,X_r^{s,X^{t,x}_s},Y_r^{s,X^{t,x}_s},Z_r^{s,X^{t,x}_s}\right)dr\right)\\
\begin{aligned}
=&\Eds\left(Y^{s,X^{t,x}_s}_s+\int_s^Tf\left(r,X_r^{s,X^{t,x}_s},Y_r^{s,X^{t,x}_s},Z_r^{s,X^{t,x}_s}\right)dr\right)\\
=&\Eds\left(Y^{t,x}_s+\int_s^Tf\left(r,X_r^{t,x},Y_r^{t,x},Z_r^{t,x}\right)dr\right)\\
=&\Eds\left(Y^{t,x}_t\right)=u(t,x),
\end{aligned}
\end{multline*}
which is just a simple consequence of the uniqueness of the solutions to the FBSDE, while for the $G$--heat equation \eqref{eq:gheat} this means that
\begin{multline}\label{eq:gheatdpp}
\sup_{\sigma\in\Acal_G}\Eds\left(u\left(s,x+\int_t^s\sigma(r)dW_r\right)\right)\\
\begin{aligned}
=&\sup_{\sigma\in\Acal_G}\sup_{\sigma'\in\Acal_G}\Eds\left(g\left(x+\int_t^s\sigma(r)dW_r+\int_s^T\sigma'(r)dW_r\right)\right)\\
=&\sup_{\sigma\in\Acal_G}\Eds\left(g\left(x+\int_t^T\sigma(r)dW_r\right)\right)\\
=&u(t,x).
\end{aligned}
\end{multline}
The proof of \eqref{eq:gheatdpp} is contained in \cite[Subsection 3.1]{art:pathpeng}. This also intuitively explain why we ask to the elements of $\Acal_F$ to be cadlag in time. We point out that in \cite{art:pathpeng} the authors ask to the elements of $\Acal_F$ to only be measurable in time, but for the analysis of the more general \cref{parprob} we will use right continuity.

The dynamic programming principle exposed in theorem \ref{cauchydpp} is a generalization of the one presented by Denis, Hu and Peng in \cite{art:pathpeng} and will be obtained in a similar way.

We start endowing $\Acal_F$ with the topology of the $L^2$--convergence on compact sets, which is to say that a sequence converges in $\Acal_F$ if and only if it converges in $L^2([0,T]\times\Omega\times K)$ for any compact $K\subset\Rds^N$, see \cref{aflp}.

Let $\varphi$ be a continuous map from $[0,T]\times\Rds^N\times\Acal_F$ into $L^2(\Omega,\Fcal_T;\Rds)$ such that, for any stopping time $\tau$ with values in $[0,T]$, $x\in\Rds^N$ and $(b,\sigma)\in\Acal_F$, $\varphi_\tau(x,b,\sigma)$ is affected only by the values that $(b,\sigma)$ assumes on $[\tau,T]$, i.e.
\[
\varphi_\tau:\Rds^N\times\Acal_F|_{[\tau,T]}\longrightarrow L^2(\Omega,\Fcal_T;\Rds),
\]
and if $(b,\sigma)\in\Acal_F^\tau$, then $\varphi_\tau(x,b,\sigma)$ is $\Fcal^\tau_T$--measurable. We further assume that, for any stopping time $\tau$ with value in $[0,T]$ and $\zeta$ in $L^2\left(\Omega,\Fcal_\tau;\Rds^N\right)$,
\begin{equation}\label{conddpp}
\sup_{(b,\sigma)\in\Acal_F}\Eds(|\varphi_\tau(\zeta,b,\sigma)|)<\infty,
\end{equation}
thus the function
\[
\Phi_\tau(\zeta):=\esss_{(b,\sigma)\in\Acal_F}\Eds(\varphi_\tau(\zeta,b,\sigma)|\Fcal_\tau)
\]
is well defined, as we will prove in \cref{peng40}.

In this section the function $\Phi$ represents, roughly speaking, the viscosity solution, $\zeta$ is the first part of a stochastic trajectory broken off at $\tau$ (this is why we restrict $\varphi_\tau$ on $\Acal_F|_{[\tau,T]}$) and $\varphi$ the function which we will use to build the viscosity solution.

\begin{lem}\label{peng40}
For each $(b_1,\sigma_1),(b_2,\sigma_2)$ in $\Acal_F$ and stopping time $\tau$ with value in $[0,T]$ there exists a $(b,\sigma)\in\Acal_F$ such that
\begin{equation}\label{eq:peng40.1}
\Eds(\varphi_\tau(\zeta,b,\sigma)|\Fcal_\tau)=\Eds(\varphi_\tau(\zeta,b_1,\sigma_1)|\Fcal_\tau)\vee\Eds(\varphi_\tau(\zeta,b_2,\sigma_2)|\Fcal_\tau).
\end{equation}
Therefore there exists a sequence $\{(b_i,\sigma_i)\}_{i\in\Nds}$ in $\Acal_F$ such that a.e.
\begin{equation}\label{eq:peng40.2}
\Eds(\varphi_\tau(\zeta,b_i,\sigma_i)|\Fcal_\tau)\uparrow\Phi_\tau(\zeta).
\end{equation}
We also have
\begin{equation}\label{eq:peng40.3}
\Eds(|\Phi_\tau(\zeta)|)\le\sup_{(b,\sigma)\in\Acal_F}\Eds(|\varphi_\tau(\zeta,b,\sigma)|)<\infty,
\end{equation}
and, for any stopping time $\tau'\le\tau$,
\begin{equation}\label{eq:peng40.4}
\Eds\left(\esss_{(b,\sigma)\in\Acal_F}\Eds(\varphi_\tau(\zeta,b,\sigma)|\Fcal_\tau)\middle|\Fcal_{\tau'}\right)=\esss_{(b,\sigma)\in\Acal_F}\Eds(\varphi_\tau(\zeta,b,\sigma)|\Fcal_{\tau'}).
\end{equation}
\end{lem}
\proof
Given $(b_1,\sigma_1),(b_2,\sigma_2)\in\Acal_F$, we define
\[
A:=\{\omega\in\Omega:\Eds(\varphi_\tau(\zeta,b_1,\sigma_1)|\Fcal_\tau)(\omega)\ge\Eds(\varphi_\tau(\zeta,b_2,\sigma_2)|\Fcal_\tau)(\omega)\},
\]
which belong to $\Fcal_\tau$, and $(b,\sigma):=(\chi_A(b_1,\sigma_1)+\chi_{A^c}(b_2,\sigma_2))\chi_{\{t\ge\tau\}}$. We thus have $(b,\sigma)\in\Acal_F$ and \eqref{eq:peng40.1}. From this, by the properties of the essential supremum, we obtain the existence of a sequence $\{(b_i,\sigma_i)\}_{i\in\Nds}$ in $\Acal_F$ such that \eqref{eq:peng40.2} is true. Furthermore, from \eqref{eq:peng40.1}, \eqref{conddpp} and \cite[Theorem 1]{art:yan}, \eqref{eq:peng40.3} and \eqref{eq:peng40.4} follow.
\endproof

\begin{rem}
To prove \cref{peng40} the randomness of the controls of $\Acal_F$ is crucial, this is the reason why we consider a set of stochastic process instead of a deterministic one.
\end{rem}

The next one is a density result on the control set $\Acal_F$.

\begin{lem}\label{densecont}
The set
\[
\Jcal^\tau:=\left\{
\begin{aligned}
&(b,\sigma)\in\Acal_F:(b,\sigma)|_{[\tau,T]}=\sum_{i=1}^n\chi_{A_i}(b_i,\sigma_i)|_{[\tau,T]},\\
&\text{ where }\{(b_i,\sigma_i)\}_{i=1}^n\subset\Acal_F^\tau\text{and }\{A_i\}_{i=0}^n\text{ is a }\Fcal_\tau\text{--partition of }\Omega
\end{aligned}
\right\}
\]
is dense in $\Acal_F$ for any stopping time $\tau$ with value in $[0,T]$.
\end{lem}
\proof
Preliminarily we will prove the \namecref{densecont} when the stopping time $\tau$ is identically equal to a $t\in[0,T]$.\\
Notice that, fixed a $k\in\Nds$, by our assumption each element of $\Acal_F$ can be approximated in $L^2([0,T]\times\Omega\times B_k(0))$ by a sequence of simple functions. We will denote with $\Bcal_k:=\Bcal([0,T]\times B_k(0))$ the Borel $\sigma$--algebra of $[0,T]\times B_k(0)$.\\
Furthermore, since the collection $\Ical$ of the rectangles $A\times B$ where $A\in\Fcal_T$ and $B\in\Bcal_k$ is a $\pi$--system which contains the complementary of its sets and generate $\sigma(\Fcal_T\times\Bcal_k)$, by \cite[Dynkin's lemma A1.3]{williams08} each set in $\sigma(\Fcal_T\times\Bcal_k)$, which is the smallest $d$--system containing $\Ical$, can be approximate by a finite union of sets in $\Ical$. Similarly, each set in $\Fcal_T$ can be approximated by finite intersection and union of sets in $\Fcal_t$ and $\Fcal_T^t$, since $\Fcal_T=\sigma\left(\Fcal_t,\Fcal_T^t\right)$.\\
Therefore, fixed $(b,\sigma)\in\Acal_F$, for any $\varepsilon>0$ there exists a simple function $s_\varepsilon$ such that $s_\varepsilon(t,\omega,x)=\sum\limits_{i=1}^n\sum\limits_{j=1}^ms^j_i(t,x)\chi_{A_i}(\omega)\chi_{A'_j}(\omega)$ where $\{A_i\}_{i=1}^n$ and $\{A'_j\}_{j=1}^m$ are respectively a $\Fcal_t$--partition and a $\Fcal^t_T$--partition of $\Omega$ and
\[
\Eds\left(\int_0^T\int_{B_k(0)}|(b,\sigma)(t,x)-s_\varepsilon(t,x)|^2dxdt\right)<\varepsilon.
\]
From this follows that for each $A_i$ and $A'_j$ there is a $\omega_i^j$, which belongs to $A_i\cap A'_j$ if $A_i\cap A'_j\ne\emptyset$, such that
\[
\sum_{i=1}^n\sum_{j=1}^m\Eds\left(\int_0^T\int_{B_k(0)}\left|(b,\sigma)\left(t,\omega_i^j,x\right)-s_i^j(t,x)\right|^2\chi_{A_i}\chi_{A_j'}dxdt\right)<\varepsilon.
\]
If we let $\left(b_i^k,\sigma_i^k\right):=\sum\limits_{j=1}^m(b,\sigma)\left(\omega_i^j\right)\chi_{A'_j}$, then $\left(b_i^k,\sigma_i^k\right)\in\Acal^t_F$ and $\left(b^k_\varepsilon,\sigma^k_\varepsilon\right):=\sum\limits_{i=1}^n\left(b^k_i,\sigma^k_i\right)\chi_{A_i}$ is an element of $\Jcal^t$ satisfying
\[
\Eds\left(\int_0^T\int_{B_k(0)}\left|(b,\sigma)(t,x)-\left(b^k_\varepsilon,\sigma^k_\varepsilon\right)(t,x)\right|^2dxdt\right)<4\varepsilon.
\]
This proves that the sequence $\left\{\left(b_{k^{-1}}^k,\sigma_{k^{-1}}^k\right)\right\}_{k\in\Nds}$ converges to $(b,\sigma)$ and thus $\Jcal^t$ is dense in $\Acal_F$.\\
Now let $\tau$ be a stopping time with values in $[0,T]$ and define
\[
\tau_n:=\left\{
\begin{aligned}
&\frac jnT,&&\text{if }\frac{j-1}nT\le\tau<\frac jnT,j\in\{1,\dotsc,n\},\\
&T,&&\text{if }\tau=T.
\end{aligned}
\right.
\]
Clearly $\tau$ is the decreasing limit of $\{\tau_n\}_{n\in\Nds}$ and for any fixed a $(b,\sigma)\in\Acal_F$ the previous step shows that there is a sequence $\left\{\left(b_i^{j,n},\sigma_i^{j,n}\right)\right\}_{i\in\Nds}\subset\Jcal^{\frac jnT}$ converging to $(b,\sigma)$ for any $n\in\Nds$ and $j\in\{1,\cdots,n\}$. Finally we define, for any $n\in\Nds$ and $(t,x)\in[0,T]\times\Rds^N$,
\begin{align*}
(b_n,\sigma_n)(t,x):=&(b,\sigma)(t,x)\chi_{\{t<\tau\}}+(b_n^{n,n},\sigma^{n,n}_n)(T,x)\chi_{\{\tau=T=t\}}\\
&+\sum_{j=1}^n\left(b_n^{j,n},\sigma^{j,n}_n\right)(\tau_n-\tau+t,x)\chi_{\left\{\frac{j-1}n\le\tau<\frac jn,t\ge\tau\right\}}.
\end{align*}
It is easy to see that $(b_n,\sigma_n)$ belongs to $\Jcal^\tau$ for any $n\in\Nds$ and it is a simple verification that $(b_n,\sigma_n)$ converges to $(b,\sigma)$. The arbitrariness of $(b,\sigma)$ hence proves that $\Jcal^\tau$ is dense in $\Acal_F$.
\endproof

\begin{lem}\label{artpeng42}
For any stopping time $\tau$ with values in $[0,T]$ and $x\in\Rds^N$, $\Phi_\tau(x)$ is deterministic. Furthermore
\begin{equation}\label{eq:artpeng42}
\Phi_\tau(x)=\esss_{(b,\sigma)\in\Acal_F}\Eds(\varphi_\tau(x,b,\sigma)|\Fcal_\tau)=\esss_{(b,\sigma)\in\Acal_F^\tau}\Eds(\varphi_\tau(x,b,\sigma)|\Fcal_\tau).
\end{equation}
\end{lem}
\proof
By \cref{densecont} $\Jcal^\tau$ is dense in $\Acal_F$, hence, by \cref{peng40}, there is a sequence
\[
\{(b_i,\sigma_i)\}_{i\in\Nds}=\left\{\sum_{j=1}^{n_i}\chi_{A_j}(b_{i,j},\sigma_{i,j})\right\}_{i\in\Nds}
\]
which belongs to $\Jcal^\tau$ such that $\Eds(\varphi_\tau(x,b_i,\sigma_i)|\Fcal_\tau)\uparrow\Phi_\tau(x)$ a.e.. But
\begin{align*}
\Eds(\varphi_\tau(x,b_i,\sigma_i)|\Fcal_\tau)\!=&\!\sum_{j=1}^{n_i}\!\chi_{A_j}\Eds(\varphi_\tau(x,b_{i,j},\sigma_{i,j})|\Fcal_\tau)\!=\!\sum_{j=1}^{n_i}\!\chi_{A_j}\Eds(\varphi_\tau(x,b_{i,j},\sigma_{i,j}))\\
\le&\max_{j=1}^{n_i}\Eds(\varphi_\tau(x,b_{i,j},\sigma_{i,j}))=\Eds(\varphi_\tau(x,b_{i,j_i},\sigma_{i,j_i})),
\end{align*}
where $j_i$ is the index that realizes the maximum in the previous formula. This implies that $\lim\limits_{i\to\infty}\Eds(\varphi_\tau(x,b_{i,j_i},\sigma_{i,j_i}))=\Phi_\tau(x)$, therefore $\Phi_\tau(x)$ is deterministic, and the inequality
\[
\lim_{i\to\infty}\Eds(\varphi_\tau(x,b_{i,j_i},\sigma_{i,j_i}))\le\esss_{(b,\sigma)\in\Acal_F^\tau}\Eds(\varphi_\tau(x,b,\sigma)|\Fcal_\tau)\le\Phi_\tau(x)
\]
confirms the equation \eqref{eq:artpeng42}.
\endproof

\begin{rem}\label{Phiisdet}
Since $\Phi_\tau(x)$ is deterministic for any stopping time $\tau$ with values in $[0,T]$ and $x\in\Rds^N$ we have, thanks to \eqref{eq:peng40.4},
\[
\Phi_\tau(x)=\sup_{(b,\sigma)\in\Acal_F}\Eds(\varphi_\tau(x,b,\sigma)).
\]
\end{rem}

\begin{lem}\label{artpeng43}
We define the function
\begin{align*}
u:[0,T]\times\Rds^N&\xrightarrow{\qquad}\Rds\\
(t,x)&\xmapsto{\qquad}\Phi_t(x)
\end{align*}
and assume that it is continuous. Then, for each stopping time $\tau$ with value in $[0,T]$ and $\zeta\in L^2\left(\Omega,\Fcal_\tau;\Rds^N\right)$, we have that $u_\tau(\zeta)=\Phi_\tau(\zeta)$ a.e..
\end{lem}
\proof
By the continuity of $u$, and consequently of $\Phi$, we only need to prove the \namecref{artpeng43} when $\tau=\sum\limits_{j=1}^n\chi_{A_j}t_j$ and $\zeta=\sum\limits_{j=1}^n\chi_{A_j}x_j$, where $t_j\in[0,T]$, $x_j\in\Eds^N$ and $\{A_j\}_{j=1}^n$ is a $\Fcal_\tau$--partition of $\Omega$. As seen in the proof of \cref{artpeng42}, for each $(t_j,x_j)$ there is a sequence $\{(b_{i,j},\sigma_{i,j})\}_{i\in\Nds}$ in $\Acal^{t_j}_F$ such that
\[
\lim_{i\to\infty}\Eds\left(\varphi_{t_j}(x_j,b_{i,j},\sigma_{i,j})\right)=\Phi_{t_j}(x_j)=u_{t_j}(x_j).
\]
Setting $(b_i,\sigma_i):=\sum\limits_{j=1}^n\chi_{A_j}(b_{i,j},\sigma_{i,j})$, we have
\begin{align*}
\Phi_\tau(\zeta)\ge&\Eds(\varphi_\tau(\zeta,b_i,\sigma_i)|\Fcal_\tau)\\
=&\sum_{j=1}^n\chi_{A_j}\Eds\left(\varphi_{t_j}(x_j,b_{i,j},\sigma_{i,j})|\Fcal_{t_j}\right)\xrightarrow[i\to\infty]{}\sum_{j=1}^n\chi_{A_j}u_{t_j}(x_j)=u_\tau(\zeta).
\end{align*}
On the other hand, for any $(b,\sigma)\in\Acal_F$,
\[
\Eds(\varphi_\tau(\zeta,b,\sigma)|\Fcal_\tau)=\sum_{j=1}^n\chi_{A_j}\Eds\left(\varphi_{t_j}(x_j,b,\sigma)|\Fcal_{t_j}\right)\le\sum_{j=1}^n\chi_{A_j}u_{t_j}(x_j)=u_\tau(\zeta),
\]
thus $\esss\limits_{b,\sigma\in\Acal_F}\Eds(\varphi_\tau(\zeta,b,\sigma)|\Fcal_\tau)\le u_\tau(\zeta)$. This completes the proof.
\endproof

\begin{rem}\label{supeqrem}
\Cref{artpeng43} says, as a consequence of \eqref{eq:artpeng42}, that
\[
\esss_{(b,\sigma)\in\Acal_F}\Eds(\varphi_\tau(\zeta,b,\sigma)|\Fcal_\tau)=\esss_{(b,\sigma)\in\Acal_F^\tau}\Eds(\varphi_\tau(\zeta,b,\sigma)|\Fcal_\tau)
\]
for any stopping time $\tau$ with value in $[0,T]$ and $\zeta\in L^2\left(\Omega,\Fcal_\tau;\Rds^N\right)$.
\end{rem}

\section{Representation Formula}\label{repform}

When $F$ is a linear operator it is known that the representation formula of its viscosity solution is built from a stochastic system, see \eqref{eq:exparsemlin}. To adapt this method to our case we will build a candidate viscosity solution using a stochastic system, defined below, and we will prove, thanks to the \cref{cauchydpp}, that it is actually a solution to \cref{parprob}.

\begin{defin}\label{fbsde}
Consider the \emph{forward backward stochastic differential equation}, FBSDE for short,
\begin{equation}\label{eq:fbsde}
\left\{\begin{aligned}
X^{t,\zeta}_s=&\zeta+\int_t^s\sigma\left(r,X^{t,\zeta}_r\right)dW_r+\int_t^sb\left(r,X^{t,\zeta}_r\right)dr,\\
Y^{t,\zeta}_s=&g\left(X^{t,\zeta}_T\right)+\int_s^Tf_{(b,\sigma)}\left(r,X_r^{t,\zeta},Y_r^{t,\zeta},Z_r^{t,\zeta}\right)dr\\
&-\int_s^TZ_r^{t,\zeta}dW_r,
\end{aligned}\right.\,s\in[t,T],
\end{equation}
where $\zeta\in L^2\left(\Omega,\Fcal_t;\Rds^N\right)$, $(b,\sigma)\in\Acal_F$ and the functions $f$ and $g$ are as in the assumptions of \cref{parprob}. Notice that the first and the second equation of \eqref{eq:fbsde} are, respectively, an SDE and a BSDE, see \cref{sdeprel} for a brief review about these topics.\\
We will call $(X,Y,Z)$ a solution to the FBSDE if $X$ is a solution to the SDE part of this system and $(Y^{t,\zeta},Z^{t,\zeta})$ is a solution to the BSDE part for any $(t,\zeta)\in[0,T]\times L^2\left(\Omega,\Fcal_t;\Rds^N\right)$.
\end{defin}

\Cref{sdeex,bsdeex} yield that, under our assumptions, there is a unique solution to \eqref{eq:fbsde}. Due to \cref{sdestopstart}, this is true even if $t$ is a stopping time. In particular from this follows the next crucial result.

\begin{prop}
For a fixed $(b,\sigma)\in\Acal_F$, let $(X,Y,Z)$ be the solution to \eqref{eq:fbsde}. Then, for any $0\le t\le r\le s\le T$,
\[
\left(X^{r,X^{t,\zeta}_r}_s,Y^{r,X^{t,\zeta}_r}_s,Z^{r,X^{t,\zeta}_r}_s\right)=\left(X^{t,\zeta}_s,Y^{t,\zeta}_s,Z^{t,\zeta}_s\right),\qquad\text{a.e.}.
\]
This holds true even if $t$, $r$ and $s$ are stopping times.
\end{prop}

For the remainder of this section, we will simply write $Y$ to denote the second term of the triplet $(X,Y,Z)$ solution to the FBSDE in \cref{fbsde}, for $(b,\sigma)$ that varies in $\Acal_F$. To ease notation we will omit the dependence of $X$, $Y$ and $Z$ on $(b,\sigma)$.

We will prove that $u(t,x):=\esss\limits_{(b,\sigma)\in\Acal_F}\Eds\left(Y^{t,x}_t\middle|\Fcal_t\right)$ is a viscosity solution to \cref{parprob}, breaking the proof in several steps. Before that however we stress out that since $Y^{t,x}_t$ is $\Fcal_t$ measurable and thanks to \cref{Phiisdet} we have that for $u$ the following identities holds true
\[
u(t,x):=\esss_{(b,\sigma)\in\Acal_F}\Eds\left(Y^{t,x}_t\middle|\Fcal_t\right)=\sup_{(b,\sigma)\in\Acal_F}\Eds\left(Y^{t,x}_t\right)=\esss_{(b,\sigma)\in\Acal_F}Y^{t,x}_t.
\]

\begin{prop}\label{eq:parsubopgenviscont}
The function $u(t,x):=\sup\limits_{(b,\sigma)\in\Acal_F}\Eds\left(Y^{t,x}_t\right)$ is $\frac12$--Hölder continuous in the first variable and Lipschitz continuous in the second one. Furthermore we have that there is a constant $c$, which depends only on $\ell$, $\mu$ and $T$, such that
\begin{equation}\label{eq:parvisbound}
\Eds\left(|u(\tau,\zeta)|^2\right)\le\sup_{(b,\sigma)\in\Acal_F}\Eds\left(\left|Y^{\tau,\zeta}_\tau\right|^2\right)\le c\left(1+\Eds\left(|\zeta|^2\right)\right),
\end{equation}
for any stopping time $\tau$ with values in $[0,T]$ and $\zeta\in L^2\left(\Omega,\Fcal_t;\Rds^N\right)$.
\end{prop}

We point out that this proposition permits us to use the results of \cref{dynprincsec} on $u$. In particular the function $(t,x,b,\sigma)\mapsto Y^{t,x}_t$ satisfies the same conditions of $\varphi$ in \cref{dynprincsec}, thanks to \cref{sdeex,bsdeex,sdemark,bsdemark}. Furthermore we prove here that it satisfies \eqref{conddpp} and the continuity of $u$, which is needed for \cref{artpeng43}. See also \cref{Phiisdet}.

\proof
Notice that Jensen's inequality yields
\begin{align*}
|u(t,x)-u(s,y)|=&\left|\sup_{(b,\sigma)\in\Acal_F}\Eds\left(Y^{t,x}_t\right)-\sup_{(b,\sigma)\in\Acal_F}\Eds\left(Y^{s,y}_s\right)\right|\\
\le&\sup_{(b,\sigma)\in\Acal_F}\left(\Eds\left(\left|Y^{t,x}_t-Y^{s,y}_s\right|^2\right)\right)^\frac12
\end{align*}
and
\[
|u(t,\zeta)|=\left|\sup_{(b,\sigma)\in\Acal_F}\Eds\left(Y^{t,\zeta}_t\right)\right|\le\left(\sup_{(b,\sigma)\in\Acal_F}\Eds\left(\left|Y^{t,\zeta}_t\right|^2\right)\right)^\frac12,
\]
for any $t,s\in[0,T]$, $x,y\in\Rds^N$ and $\zeta\in L^2\left(\Omega,\Fcal_t;\Rds^N\right)$.\\
In order to prove that $u$ is Hölder continuous in the first variable we assume that $s\le t$, hence we obtain from \cref{bsdeex,sdeex} that there exist two constants $c_1$ and $c_2$, which depend on $\ell$, $\mu$ and $T$, such that
\begin{multline*}
\mathds E\left(\left|Y^{t,x}_t-Y^{s,x}_s\right|^2\right)\le2\mathds E\left(\left|Y^{s,x}_t-Y^{s,x}_s\right|^2\right)+2\mathds E\left(\left|Y^{t,x}_t-Y^{s,x}_t\right|^2\right)\\
\begin{aligned}
=&2\mathds E\left(\left|Y^{s,x}_t-Y^{s,x}_s\right|^2\right)+2\mathds E\left(\left|Y^{t,x}_t-Y^{t,X^{s,x}_t}_t\right|^2\right)\\
\le&c_1\mathds E\biggl(
\begin{aligned}[t]
&\int_s^t(|f_{(b,\sigma)}\left(r,X^{s,x}_r,Y^{s,x}_r,Z^{s,x}_r)|^2+|Z_r^{s,x}|^2\right)dr\\
&+\left|g\left(X^{t,x}_T\right)-g\left(X^{t,X^{s,x}_t}_T\right)\right|^2\\
&+\left.\int_t^T\left|f_{(b,\sigma)}\left(r,X^{t,x}_r,Y^{t,x}_r,Z^{t,x}_r\right)-f_{(b,\sigma)}\left(r,X^{t,X^{s,x}_t}_r,Y^{t,x}_r,Z^{t,x}_r\right)\right|^2dr\right)
\end{aligned}\\
\le&c_2\mathds E\biggl(
\begin{aligned}[t]
&\int_s^t\left(1+\left|X^{s,x}_r\right|^2+\left|Y^{s,x}_r\right|^2+\left|Z^{s,x}_r\right|^2\right)dr\\
&+\left.\left|X^{t,x}_T-X^{t,X^{s,x}_t}_T\right|^2+\int_t^T\left|X^{t,x}_r-X^{t,X^{s,x}_t}_r\right|^2dr\right).
\end{aligned}
\end{aligned}
\end{multline*}
Using once again \cref{bsdeex,sdeex}, we have that there exist five constants $c_3$, $c_4$, $c_5$, $c_6$ and $c_7$ depending upon $\mu$, $\ell$ and $T$ such that
\begin{equation}\label{eq:auxbondparsolaux}
\begin{aligned}
\mathds E\left(\int_s^t\left(|Y^{s,x}_r|^2+|Z^{s,x}_r|^2\right)dr\right)\le&c_3\mathds E\left(\int_s^t\left|f_{(b,\sigma)}(r,X^{s,x}_r,0,0)\right|^2dr\right)\\
\le&c_4\mathds E\left(\int_s^t\left(1+|X^{s,x}_r|^2\right)dr\right)\\
\le&c_5(t-s)\left(1+|x|^2\right),
\end{aligned}
\end{equation}
and
\[
\mathds E\left(\sup_{r\in[t,T]}\left|X^{t,x}_r-X^{t,X^{s,x}_t}_r\right|^2\right)\le c_6\mathds E\left(|x-X^{s,x}_t|^2\right)\le c_7(t-s).
\]
These inequalities together imply that
\[
\mathds E\left(\left|Y^{t,x}_t-Y^{s,x}_s\right|^2\right)\le c_8(t-s)\left(1+|x|^2\right),
\]
where $c_8$ is a constants depending only on $\mu$, $\ell$ and $T$, showing that $u$ is $\frac12$--Hölder continuous in the first variable. The rest can be proved similarly.
\endproof

We can now prove the dynamic programming principle for $u$.

\begin{theo}[Dynamic programming principle]\label[dpp]{cauchydpp}
For any $(b,\sigma)\in\Acal_F$ we let $\left(\overline Y,\overline Z\right)$ be the solution to the BSDE
\begin{equation}\label{eq:cauchydpp1}
\overline Y_s=u\left(\tau,X_\tau^{t,x}\right)+\int_{s\wedge\tau}^\tau f_{(b,\sigma)}\left(r,X^{t,x}_r,\overline Y_r,\overline Z_r\right)dr-\int_{s\wedge\tau}^\tau\overline Z_rdW_r,
\end{equation}
where $s\in[t,T]$ and $\tau$ is a stopping time with values in $[t,T]$. Then we have $\sup\limits_{(b,\sigma)\in\Acal_F}\Eds\left(\overline Y_t\right)=u(t,x)$.
\end{theo}
\proof
Fixed $(\overline b,\overline\sigma)\in\Acal_F$ in \eqref{eq:cauchydpp1} we define $\overline X:=X_{(b,\sigma)}$ and the subset of $\Acal_F$
\[
\overline{\Acal}_F:=\left\{(b,\sigma)\in\Acal_F:(b,\sigma)(s)=\left(\overline b,\overline\sigma\right)(s)\text{ for any }s\in[t,\tau)\right\}.
\]
From \cref{artpeng43} we know that
\[
\esss_{(b,\sigma)\in\overline{\Acal}_F}\Eds\left(Y_\tau^{t,x}\middle|\Fcal_\tau\right)=\esss_{(b,\sigma)\in\overline{\Acal}_F}\Eds\left(Y_\tau^{\tau,\overline X^{t,x}_\tau}\middle|\Fcal_\tau\right)=u\left(\tau,\overline X_\tau^{t,x}\right)
\]
and \cref{peng40} yields the existence of a sequence $\{(b_n,\sigma_n)\}_{n\in\Nds}$ in $\overline{\Acal}_F$ and a corresponding sequence $\{Y_n\}_{n\in\Nds}$ such that
\[
\lim_{n\to\infty}\Eds\left(Y_{n,\tau}^{t,x}\middle|\Fcal_\tau\right)=\esss_{(b,\sigma)\in\overline{\Acal}_F}\Eds\left(Y_\tau^{t,x}\middle|\Fcal_\tau\right)=u\left(\tau,\overline X_\tau^{t,x}\right).
\]
Then, by \cref{bsdeex} and the dominated convergence theorem, there exists a constant $c$ such that
\[
\lim_{n\to\infty}\Eds\left(\left|\overline Y_t-Y^{t,x}_{n,t}\right|^2\right)\le\lim_{n\to\infty}c\Eds\left(\left|u\left(\tau,\overline X_\tau^{t,x}\right)-Y^{t,x}_{n,\tau}\right|^2\right)=0,
\]
hence, up to subsequences,
\begin{equation}\label{eq:cauchydpp2}
\lim\limits_{n\to\infty}\Eds\left(Y_{n,t}^{t,x}\right)=\Eds\left(\overline Y_t\right).
\end{equation}
Furthermore, thanks to \cref{bsdeconf}, $Y^{t,x}_t\le\overline Y_t$ for any $(b,\sigma)\in\overline{\Acal}_F$, which together with \eqref{eq:cauchydpp2} implies that $\sup\limits_{(b,\sigma)\in\overline{\Acal}_F}\Eds\left(Y^{t,x}_t\right)=\Eds\left(\overline Y_t\right)$. Therefore we can use the arbitrariness of $\left(\overline b,\overline\sigma\right)$ to obtain our conclusion:
\[
\sup_{(b,\sigma)\in\Acal_F}\Eds\left(\overline Y_t\right)=\sup_{(b,\sigma)\in\Acal_F}\Eds\left(Y_t^{t,x}\right)=u(t,x).
\]
\endproof

We conclude this \namecref{repform} with our main statement.

\begin{theo}
The function $u(t,x):=\sup\limits_{(b,\sigma)\in\Acal_F}\Eds\left(Y^{t,x}_t\right)$ is the only viscosity solution to \cref{parprob} satisfying polynomial growth condition such that $u(T,x)=g(x)$ for any $x$ in $\Rds^N$.
\end{theo}
\proof
By \eqref{eq:parvisbound} $u$ has polynomial growth, thus the uniqueness of $u$ follows from \cref{compparprob}. Moreover, from \cref{eq:parsubopgenviscont}, we know that $u$ is continuous and it is easy to see that $u(T,x)=g(x)$ for any $x\in\Rds^N$, so we only need to prove that it is a viscosity solution.\\
We start showing that $u$ is a subsolution. In order to do so we preliminarily fix $(t,x)\in(0,T)\times\Rds^N$ and $(b,\sigma)\in\Acal_F^t$, and observe that $(b,\sigma)(t)$ is $\Fcal_t^t$--measurable. By Blumenthal 0--1 law then it is a.e. deterministic, which is to say that there is a $(b',\sigma')\in K_F$ such that $(b,\sigma)(t,y)=(b',\sigma')(t,y)$ a.e. for any $y\in\Rds^N$. Similarly for the corresponding $Y$ we have $Y^{t,y}_t=\Eds\left(Y^{t,y}_t\right)$ a.e. for any $y\in\Rds^N$. We will prove that $v(s,y):=\Eds\left(Y^{s,y}_s\right)$ is a subsolution at $(t,x)$ and, since by \cref{supeqrem,Phiisdet} $u(t,x)=\sup\limits_{(b,\sigma)\in\Acal_F^t}\Eds\left(Y^{t,x}_t\right)$, well known properties of viscosity solution yield that $u$ is a subsolution at $(t,x)$. From the arbitrariness of $(t,x)$ will then follows that $u$ is a subsolution to \cref{parprob}.\\
Let $\varphi$ be a supertangent to $v$ at $(t,x)$ which we assume, without loss of generality, equal to $v$ at $(t,x)$ and, by contradiction, such that
\begin{multline*}
\partial_t\varphi(t,x)+\frac12\left\langle\sigma\sigma^\dag,D^2_x\varphi\right\rangle(t,x)+(\nabla_x\varphi b)(t,x)+f_{(b,\sigma)}(t,x,v,\nabla_x\varphi\sigma)\\
\le\partial_t\varphi(t,x)+F\left(t,x,v,\nabla_x\varphi,D^2_x\varphi\right)<0.
\end{multline*}
We define the stopping time
\begin{equation}\label{eq:solprob1}
\tau:=\inf\left\{s\in[t,T]:Y^{t,x}_s=Y^{s,X^{t,x}_s}_s>\varphi\left(s,X_s^{t,x}\right)\right\}
\end{equation}
and observe that $\tau>t$ on a set of positive measure, since otherwise we would have that $\varphi\left(s,X^{t,x}_s\right)<v\left(s,X^{t,x}_s\right)$ for any $s>t$ in some neighborhood of $t$. Since the controls in $\Acal_F$ are right continuous we also assume, eventually taking a smaller stopping time $\tau$, that, for any $s\in[t,T]$,
\begin{multline}\label{eq:solprob2}
\partial_t\varphi\left(s\wedge\tau,X^{t,x}_{s\wedge\tau}\right)+\frac12\left\langle\sigma\sigma^\dag,D^2_x\varphi\right\rangle\left(s\wedge\tau,X^{t,x}_{s\wedge\tau}\right)\\
+(\nabla_x\varphi b)\left(s\wedge\tau,X^{t,x}_{s\wedge\tau}\right)+f_{(b,\sigma)}\left(s\wedge\tau,X^{t,x}_{s\wedge\tau},v,\nabla_x\varphi\sigma\right)<0.
\end{multline}
Now define $\left(\overline Y_s,\overline Z_s\right):=\left(Y^{t,x}_{s\wedge\tau},Z^{t,x}_{s\wedge\tau}\right)$, which solves the BSDE
\[
\overline Y_s=Y_\tau^{t,x}+\int_{s\wedge\tau}^\tau f_{(b,\sigma)}\left(r,X^{t,x}_r,\overline Y_r,\overline Z_r\right)dr-\int_{s\wedge\tau}^\tau\overline Z_rdW_r,\quad s\in[t,T],
\]
and $\left(\hat Y_s,\hat Z_s\right):=\left(\varphi\left(s,X^{t,x}_{s\wedge\tau}\right),(\nabla_x\varphi\sigma)\left(s,X^{t,x}_{s\wedge\tau}\right)\right)$, which by Itô's formula is solution, for any $s\in[t,T]$, to
\begin{align*}
\hat Y_s=&\varphi\left(\tau,X^{t,x}_\tau\right)-\int_{s\wedge\tau}^\tau\hat Z_rdW_r\\
&-\int_{s\wedge\tau}^\tau\!\left(\partial_t\varphi\left(r,X^{t,x}_r\right)+\frac12\left\langle\sigma\sigma^\dag,D^2_x\varphi\right\rangle\left(r,X^{t,x}_r\right)+(\nabla_x\varphi b)\left(r,X^{t,x}_r\right)\right)dr.
\end{align*}
Thanks to \eqref{eq:solprob1} and \eqref{eq:solprob2}, we obtain from \cref{bsdeconf,stopbsde} that $\varphi(t,x)>v(t,x)$, in contradiction to our assumptions. Therefore $u$ is a subsolution to \cref{parprob}.\\
We will now prove that $u$ is also a supersolution. Fixed a $(t,x)$ in $(0,T)\times\Rds^N$, let $\psi$ be a subtangent to $u$ in $(t,x)$ such that $\psi(t,x)=u(t,x)$. We know, thanks to \cref{probass2}, that there is a continuous and deterministic $(b,\sigma)\in\Acal_F$ for which
\begin{align*}
F\left(t,x,u,\nabla_x\psi,D^2_x\psi\right)=&\frac12\left\langle\sigma\sigma^\dag,D^2_x\psi\right\rangle(t,x)+(\nabla_x\psi b)(t,x)\\
&+f_{(b,\sigma)}(t,x,u,\nabla_x\psi\sigma),
\end{align*}
thus we assume by contradiction
\begin{multline*}
\partial_t\psi(t,x)+\frac12\left\langle\sigma\sigma^\dag,D^2_x\psi\right\rangle(t,x)+(\nabla_x\psi b)(t,x)+f_{(b,\sigma)}(t,x,u,\nabla_x\psi\sigma)\\
=\partial_t\psi(t,x)+F\left(t,x,u,\nabla_x\psi,D^2_x\psi\right)>0.
\end{multline*}
We then have by our assumptions that there is a $\delta>0$ such that
\begin{equation}\label{eq:solprob3}
\partial_t\psi(s,y)+\frac12\left\langle\sigma\sigma^\dag,D^2_x\psi\right\rangle(s,y)+(\nabla_x\psi b)(s,y)+f_{(b,\sigma)}(s,y,u,\nabla_x\psi\sigma)\!<\!0
\end{equation}
and
\begin{equation}\label{eq:solprob4}
\psi(s,y)\le u(s,y)
\end{equation}
for any $(s,y)\in[t,t+\delta)\times B_\delta(x)$, therefore the stopping time
\[
\tau:=(t+\delta)\wedge\inf\left\{s\in[t,T]:\left|X_s^{t,x}-x\right|\ge\delta\right\}
\]
is a.e. bigger than $t$. Let $X$ be the solution to the SDE $(b,\sigma)$, $\left(\overline Y_s,\overline Z_s\right)$ be the solution to the BSDE
\[
\overline Y_s=u\left(\tau,X_\tau^{t,x}\right)+\int_{s\wedge\tau}^\tau f_{(b,\sigma)}\left(r,X^{t,x}_r,\overline Y_r,\overline Z_r\right)dr-\int_{s\wedge\tau}^\tau\overline Z_rdW_r,\quad s\in[t,T]
\]
and $\left(\hat Y_s,\hat Z_s\right):=\left(\psi\left(s,X^{t,x}_{s\wedge\tau}\right),(\nabla_x\psi\sigma)\left(s,X^{t,x}_{s\wedge\tau}\right)\right)$ which, by Itô's formula, is solution, for any $[s\in[t,T]$, to
\[
\begin{aligned}
\hat Y_s=&\psi\left(\tau,X^{t,x}_\tau\right)-\int_{s\wedge\tau}^\tau\hat Z_rdW_r\\
&-\int_{s\wedge\tau}^\tau\!\left(\partial_t\psi\left(r,X^{t,x}_r\right)+\frac12\left\langle\sigma\sigma^\dag,D^2_x\psi\right\rangle\left(r,X^{t,x}_r\right)+(\nabla_x\psi b)\left(r,X^{t,x}_r\right)\right)\!dr.
\end{aligned}
\]
The \cref{cauchydpp} finally yields
\begin{equation}\label{eq:eqassum}
\sup_{(b,\sigma)\in\Acal_F}\Eds\left(\overline Y_t\right)=u(t,x)=\psi(t,x),
\end{equation}
which is in contradiction with \eqref{eq:solprob3}, since together with \eqref{eq:solprob4} it implies, thanks to \cref{bsdeconf,stopbsde}, that $\overline Y_t>\psi(t,x)$ a.e.. Then the arbitrariness of $(t,x)$ proves that $u$ is a supersolution to \cref{parprob} and concludes the proof.
\endproof

\begin{appendices}

\section{Connection with 2BSDEs}\label{2bsdesec}

In \cite{int2bsde} it is established a connection between second order BSDEs, 2BSDEs for short, and parabolic PDEs similar to the ones studied here, which can be used to obtain representation formulas. The link between this approach and this method can be seen, in particular, in the dynamic programming principle we have developed. We emphasize, on the other side, starting from the dynamic programming principle, we can provide a new angle for looking at 2BSDE. We present this point of view in the \namecref{2bsdesec}. This is clearly just a short note and not a complete analysis.

We start giving the formulation of 2BSDE as in \cite{sonertouzizhang11,possamai13}. Assume that $\Omega:=\left\{\omega\in C\left([0,T];\Rds^N\right):\omega_0=0\right\}$ and let $\Pds_0$ be the Wiener measure. Note that in this space the Brownian motion $W$ is a projection, i.e. $W_t(\omega)=\omega_t$. Denote with $[W]_t$ the quadratic variation of the projection and with
\[
\hat a_t:=\lim_{\varepsilon\downarrow0}\frac{[W]_t-[W]_{t-\varepsilon}}{\varepsilon}
\]
its variation. We will then denote with $\Pcal_W$ the set of the probability measures $\Pds$ such that $[W]$ is absolutely continuous in $t$ and $\hat a\in\Sds^N_+$, $\Pds$--a.e.. In particular $\Pds_0\in\Pcal_W$ because $[W]_t=tI_N$ and $\hat a_t=I_N$, $\Pds_0$--a.e., where $I_N$ is the $N\times N$ identity matrix. Moreover, let $\Pcal_S$ be the subset of $\Pcal_W$ composed by the probability measures $\Pds^\alpha:=\Pds_0\circ(X^\alpha)^{-1}$, where
\[
X^\alpha_t:=\int_0^t\alpha^{1/2}_sdW_s,\qquad\Pds_0\text{--a.e.},
\]
and $\alpha$ is a progressive process in $\Sds^N_+$ such that, for two fixed $\underline a,\overline a$ in $\Sds^N_+$, $\underline a\le\alpha\le\overline a$, $\Pds_0$--a.e.. It is then apparent a link between $\Pcal_S$ and the control set $\Acal_F$.\\
Now, given a function
\[
h:[0,T]\times\Omega\times\Rds\times\Rds^N\times D_h\to\Rds,
\]
where $D_h$ is a subset of $\Rds^{N\times N}$ containing 0, define for any $a\in\Sds^N_+$
\[
f(t,\omega,y,z,a):=\sup_{\gamma\in D_h}\left(\frac12\langle\gamma,a\rangle-h(t,\omega,y,z,\gamma)\right).
\]
Furthermore let $\hat f(t,y,z):=f(t,y,z,\hat a_t)$,
\[
\Pcal_h^2:=\left\{\Pds\in\Pcal_S:\Eds^{\Pds}\left(\int_0^T\left|\hat f(t,0,0)\right|^2dt\right)<\infty\right\}
\]
and $\Pcal_h^2(t,\Pds):=\left\{\Pds'\in\Pcal_h^2:\Pds'=\Pds\text{ on }\Fcal_t\right\}$.

A pair of progressive processes $(Y,Z)$ is solution to the 2BSDE
\begin{equation}\label{eq:2bsde}
Y_t=\xi+\int_t^T\hat f(s,Y_s,Z_s)ds-\int^T_tZ_sdW_s+K_T-K_t
\end{equation}
if, for any $\Pds\in\Pcal_h^2$,
\begin{myenum}
\item $Y_T=\xi$, $\Pds$--a.e.;
\item the process $K^{\Pds}$ defined below has non decreasing path $\Pds$--a.e.,
\[
K_t^{\Pds}=Y_0-Y_t-\int_0^t\hat f(s,Y_s,Z_s)ds+\int^t_0Z_sdW_s,\quad t\in[0,T],\,\Pds\text{--a.e.};
\]
\item the family $\left\{K^{\Pds},\Pds\in\Pcal^2_h\right\}$ satisfies the minimum condition
\begin{equation}\label{eq:kmincond}
K_t^{\Pds}=\essi_{\Pds'\in\Pcal_h^2(t,\Pds)}\Eds^{\Pds}\left(K_T^{\Pds'}\middle|\Fcal_t\right),\qquad\Pds\text{--a.e.}.
\end{equation}
\end{myenum}
Under suitable conditions the 2BSDE \eqref{eq:2bsde} admits a unique solution and, if we denote with $\left(Y^{\Pds}(r,\xi),Z^{\Pds}(r,\xi)\right)$ the solution to the BSDE
\[
Y_t=\xi+\int_t^r\hat f(s,Y_s,Z_s)ds-\int^r_tZ_sdW_s,\qquad t\in[0,r],\,\Pds\text{--a.e.},
\]
it can be proved that $Y_t=\esss\limits_{\Pds'\in\Pcal_h^2(t,\Pds)}Y_t^{\Pds'}(r,Y_r)$ for any $r\in[t,T]$ and $\Pds\in\Pcal_h^2$. The last identity is a dynamic programming principle and can be seen as the connection between 2BSDE and our method.

Now we will show a different formulation of 2BSDEs, using controls instead of probability measures. Let $\Acal$ be a control set, made up by the progressive processes in $L^2([0,T]\times\Omega,\mathds P_0;B)$, where $B$ is a Banach space, and, for any $t\in[0,T]$,
\[
\Acal(t,\alpha):=\left\{\alpha'\in\Acal:\alpha'_s=\alpha_s\text{ for any }s\in[0,t]\right\}.
\]
Then define the function
\[
f:[0,T]\times\Omega\times\Rds\times\Rds^N\times B\to\Rds
\]
and assume that there exists a $C>0$ such that
\[
|f(t,y,z,\alpha)-f(t,y,z,\alpha')|\le C|\alpha-\alpha'|.
\]
We will also assume that, for each $\alpha\in\Acal$, $f_\alpha(t,y,z):=f(t,y,z,\alpha_t)$ satisfies \cref{bsdeassum} uniformly with respect to $\alpha$. We point out that these conditions are not intended to be minimal. We then have that the BSDEs
\[
Y^\alpha_t=\xi_\alpha+\int_t^Tf_\alpha(s,Y_s^\alpha,Z_s^\alpha)ds-\int_t^TZ^\alpha_sdW_s,\qquad t\in[0,T],
\]
admit a unique solution for any $\alpha\in\Acal$. If we moreover require that, for any $t\in[0,T]$, $\sup\limits_{\alpha\in\Acal}\Eds\left(|Y^\alpha_t|^2\right)<\infty$ (which can be achieved if, for example, $B$ is compact), then reasoning as in the proof of \cref{cauchydpp},
\begin{equation}\label{eq:poordpp}
\overline Y_t^\alpha:=\esss_{\alpha'\in\Acal(t,\alpha)}Y^{\alpha'}_t=\esss_{\alpha'\in\Acal(t,\alpha)}Y^{\alpha'}_t\left(r,\overline Y_r^{\alpha'}\right),\quad\text{for any }0\le t\le r\le T,
\end{equation}
where $Y^{\alpha'}\left(r,\overline Y_r^{\alpha'}\right)$ is solution to the BSDE $\left(\overline Y^{\alpha'}_r,f_\alpha,r\right)$. It is easy to see that, for any $\alpha\in\Acal$, $\overline Y^\alpha$ is a continuous progressive process in $L^2$ and $\overline Y^\alpha_T=\xi_\alpha$ a.e..\\
Using the same arguments of \cite{possamai13} we have that, for each $\alpha\in\Acal$, there exist two progressive processes in $L^2$, $Z^\alpha$ and $K^\alpha$, such that $K^\alpha$ is a continuous and increasing process in $t$ with $K_0^\alpha=0$ and
\[
\overline Y_t^\alpha=\xi_\alpha+\int_t^Tf_\alpha\left(s,\overline Y^\alpha_s,\overline Z_s^\alpha\right)ds-\int^T_t\overline Z_s^\alpha dW_s+K^\alpha_T-K^\alpha_t,\qquad t\in[0,T].
\]
We also have that, as in \eqref{eq:kmincond},
\[
K_t^\alpha=\essi_{\alpha'\in\Acal(t,\alpha)}\Eds\left(K_T^{\alpha'}\middle|\Fcal_t\right),\qquad\text{a.e. for any }t\in[0,T],\,\alpha\in\Acal.
\]

If we let $X$ be as in \eqref{eq:fbsde} and $\left(\overline Y,\overline Z,K\right)$ be the solution to the 2BSDE (we omit the dependence on the control set $\Acal_F$)
\[
\overline Y_s^{t,x}=g\left(X^{t,x}_T\right)+\int_s^Tf_{(b,\sigma)}\left(r,X^{t,x}_T,\overline Y_r^{t,x},\overline Z_r^{t,x}\right)ds-\int^T_s\overline Z_r^{t,x}dW_r+K^{t,x}_T-K^{t,x}_s
\]
for any $x\in\Rds^N$ and $0\le t\le s\le T$, we then have that the function $u(t,x):=\overline Y^{t,x}_t$ is the viscosity solution to \cref{parprob}.

\section{Stochastic Differential Equations}\label{sdeprel}

In this \namecref{sdeprel} are given some results on \emph{stochastic differential equations}, SDEs for short, we use in this paper.\\
Consider the SDE $(b,\sigma)$
\begin{equation}\label{eq:sde}
X^{t,\zeta}_s=\zeta+\int_t^s\sigma\left(r,X^{t,\zeta}_r\right)dW_r+\int_t^sb\left(r,X^{t,\zeta}_r\right)dr,\qquad s\in[t,\infty),
\end{equation}
under the following \namecref{sdeassum}:

\begin{assum}\label{sdeassum}
Let $t\in[0,\infty)$, $\zeta\in L^2\left(\Omega,\Fcal_t;\Rds^N\right)$ and
\[
b:[0,\infty)\times\Omega\times\Rds^N\to\Rds^N\text{ and }\sigma:[0,\infty)\times \Omega\times\Rds^N\to\Rds^{N\times M}.
\]
Assume there exists a positive constant $\ell$ such that
\begin{myenum}[series=sdeassum]
\item $\{(b,\sigma)(t,x)\}_{t\in[0,\infty)}$ is a progressive process belonging to $L^2([0,T]\times\Omega)$ for any $x\in\Rds^N$ and $T\in[0,\infty)$
\end{myenum}
and a.e., for any $r\in[0,\infty)$, $x,x'\in\Rds^N$,
\begin{myenum}[resume=sdeassum]
\item $|b(r,x)-b(r,x')|+|\sigma(r,x)-\sigma(r,x')|\le\ell|x-x'|$.
\end{myenum}
\end{assum}

A solution to this SDE is a continuous progressive process $X$ as in \eqref{eq:sde} such that $X\in L^2([0,T]\times\Omega)$ for any $T\in[0,\infty)$. The next theorem summarizes some SDE result given in \cite{krylovbook}.

\begin{theo}\label{sdeex}
Under \cref{sdeassum} there exists a unique solution to the SDE \eqref{eq:sde}. Moreover, if $\overline X$ is solution to the SDE $\left(\overline b,\overline\sigma\right)$, there is for each $T\in[t,\infty)$ a constant $c$, depending only on $\ell$ and $T$ such that
\begin{multline*}
\Eds\left(\sup_{s\in[t,T]}\left|X^{t,\zeta}_s-\overline X^{t,\zeta'}_s\right|^2\right)\!\le\!c\Eds\!\left(|\zeta-\zeta'|^2\!+\!\int_t^T\!\left|b\!\left(s,X_s^{t,\zeta}\right)\!-\overline b\left(s,X_s^{t,\zeta}\right)\right|^2\!ds\!\right)\\
+c\Eds\left(\int_t^T\left|\sigma\left(s,X_s^{t,\zeta}\right)-\overline\sigma\left(s,X_s^{t,\zeta}\right)\right|^2ds\right).
\end{multline*}
\end{theo}

\begin{rem}\label{sdestopstart}
The results obtained in this section hold even for SDEs with an a.e. finite stopping time $\tau$ as starting time. Indeed if for any $\zeta$ in $L^2\left(\Omega,\Fcal_\tau;\Rds^N\right)$ we define
\[
\overline b(t,x):=b(t,x+\zeta)\chi_{\{\tau\le t\}}\text{ and }\,\overline\sigma(t,x):=\sigma(t,x+\zeta)\chi_{\{\tau\le t\}},
\]
then $X^{\tau,\zeta}$ is solution of the SDE $(b,\sigma)$ if and only if $\overline X^{0,0}:=X^{\tau,\zeta}-\zeta$ is solution of the SDE $\left(\overline b,\overline\sigma\right)$. The claim can be easily obtained from this.
\end{rem}

\begin{rem}\label{sdemark}
By the strong Markov property, for any a.e. finite stopping time $\tau$, the process $\{W_t^\tau\}_{t\in[0,\infty)}:=\{W_{\tau+t}-W_\tau\}_{t\in[0,\infty)}$ is a Brownian motion. Thus if $b$ and $\sigma$ are are progressive with respect to the filtration $\{\Fcal^\tau_t\}_{t\in[0,\infty)}$ then any solution to the SDE $(b,\sigma)$ with initial data $\tau+t$ and $\zeta\in L^2\left(\Omega,\Fcal_t^\tau;\Rds^N\right)$ is also progressive with respect to that filtration. In fact, in this case, the stochastic integral with respect to $W_t^\tau$ is the same as the one with respect to $W_{\tau+t}$.
\end{rem}

\subsection{Backward Stochastic Differential Equations}

Here we give some results on \emph{backward stochastic differential equations}, BSDEs for short, used in our investigation. Most of them are well known and actually hold under more general assumptions. We refer to \cite{syspardoux,phamsdebook,pardouxbook,brianddelyonhupardouxstoica03} for a detailed overview on this topic.

We will work under the followings assumptions:

\begin{assum}\label{bsdeassum}
Let $T\in[0,\infty)$, $\xi\in L^2\left(\Omega,\Fcal_T,\Pds;\Rds\right)$ and
\[
f:[0,T]\times \Omega\times\Rds\times\Rds^N\to\Rds.
\]
Assume there exists a positive constant $\ell$ and a real number $\mu$ such that
\begin{myenum}[series=bsdeassum]
\item $\{f(s,y,z)\}_{s\in[0,T]}$ is a progressive process belonging to $L^2([0,T]\times\Omega)$ for any $y\in\Rds$ and $z\in\Rds^N$
\end{myenum}
and a.e., for any $t\in[0,T]$, $y,y'\in\Rds$ and $z,z'\in\Rds^N$,
\begin{myenum}[resume=bsdeassum]
\item $|f(t,y,z)|\le|f(t,0,0)|+\ell(1+|y|+|z|)$;
\item $|f(t,y,z)-f(t,y,z')|\le\ell|z-z'|$;
\item $(y-y')(f(t,y,z)-f(t,y',z))\le\mu|y-y'|^2$;
\item $v\mapsto f(t,v,z)$ is continuous.
\end{myenum}
\end{assum}

A solution to the BSDE $(\xi,f,T)$, where $\xi$ and $T$ have respectively the role of a final condition and a terminal time, is a pair $(Y,Z)$ of continuous progressive processes belonging to $L^2([0,T]\times\Omega)$ such that
\begin{equation}\label{eq:bsde}
Y_t=\xi+\int_t^Tf(s,Y_s,Z_s)ds-\int_t^TZ_sdW_s,\quad\text{for any }t\in[0,T].
\end{equation}

The followings are classical results of BSDE theory, see for example theorems 1.2 and 1.3, and proposition 1.3 in \cite{syspardoux}.

\begin{theo}\label{bsdeex}
Under the \cref{bsdeassum} the BSDE \eqref{eq:bsde} has a unique solution $(Y,Z)$. Furthermore, if $(Y',Z')$ is the solution to the BSDE $(\xi',f',T)$, there exists a constant $c$, which depends on $T$, $\mu$ and $\ell$, such that
\begin{multline*}
\Eds\left(\sup_{t\in[0,T]}|Y_t-Y'_t|^2+\int_0^T|Z_t-Z'_t|^2dt\right)\\
\le c\Eds\left(|\xi-\xi'|^2+\int_0^T|f(t,Y'_t,Z'_t)-f'(t,Y'_t,Z'_t)|^2dt\right).
\end{multline*}
\end{theo}

\begin{rem}\label{bsdemark}
Also for BSDEs holds a measurability result similar to the one seen in \cref{sdemark}. Indeed, if $\tau$ is a stopping time with values in $[0,T]$, $\xi\in L^2\left(\Omega,\Fcal_T^\tau,\Pds;\Rds^M\right)$ and $f$ is progressive with respect to the filtration $\{\Fcal^\tau_t\}_{t\in[0,\infty)}$, the solution $(Y,Z)$ to the BSDE $(\xi,f,T)$ restricted to $[\tau,T]$ is progressive with respect to the filtration $\{\Fcal^\tau_t\}_{t\in[0,\infty)}$.
\end{rem}

\begin{theo}\label{bsdeconf}
Let $(Y,Z)$ be the solution to the BSDE $(\xi,f,T)$ under the \cref{bsdeassum} and
\[
Y'_t=\xi'+\int_t^TV_sds-\int_t^TZ'_sdW_s,\qquad t\in[0,T],
\]
where $\xi'\in L^2\left(\Omega,\Fcal_T,\Pds;\Rds\right)$, $Y',V\in L^2([0,T]\times\Omega)$ and $Z'\in L^2([0,T]\times\Omega)$. Suppose that $\xi\le\xi'$ a.e. and $f(t,Y'_t,Z'_t)\le V_t$ a.e. for the $dt\times d\Pds$ measure. Then, for any $t\in[0,T]$, $Y_t\le Y'_t$ a.e..\\
Furthermore, if $f$ is right continuous on $[0,T]$ and $Y_0=Y'_0$ a.e., then $Y_t=Y'_t$ a.e. for any $t\in[0,T]$. In particular, whenever either $\Pds(\{\xi<\xi'\})>0$ or $f(s,Y'_s,Z'_s)<V_s$, for any $(y,z)\in\Rds\times\Rds^N$ on a set of positive $dt\times d\Pds$ measure, then $Y_0<Y'_0$.
\end{theo}

\begin{prop}\label{stopbsde}
Let $(Y,Z)$ be the solution to the BSDE \eqref{eq:bsde} under the \cref{bsdeassum} and assume that there exists a stopping time $\tau$ such that $\tau\le T$, $\xi$ is $\Fcal_\tau$--measurable and $f(t,y,z)=0$ on the set $\{\tau\le t\}$. Then $Y_t=Y_{\tau\wedge t}$ a.e. and $Z_t=0$ a.e. on the set $\{\tau\le t\}$. In particular
\[
Y_{\tau\wedge t}=\xi+\int_{\tau\wedge t}^\tau f(s,Y_s,Z_s)ds-\int_{\tau\wedge t}^\tau Z_sdW_s,\qquad\text{for any }t\in[0,T].
\]
\end{prop}

\section{Comparison Theorem}\label{compsec}

Consider the parabolic problem
\begin{equation}\label{eq:compparpde}
\partial_tu(t,x)+G\left(t,x,u,\nabla_xu,D_x^2u\right)=0,\quad t\in(0,T),x\in\Rds^N,
\end{equation}
under the following \namecref{compass}.

\begin{assum}\label{compass}
$G$ is a continuous elliptic operator which admits, for any $t\in[0,T]$, $x,p,p'\in\Rds^N$, $r,r'\in\Rds$ and $S,S'\in\Sds^N$, a positive $\ell$ and a $\mu\in\Rds$ such that
\begin{myenum}
\item\label{compass.1} $|G(t,x,r,p,S)-G(t,x,r,p,S')|\le\ell\left(1+|x|^2\right)|S-S'|$;
\item\label{compass.2} $|G(t,x,r,p,S)-G(t,x,r,p',S)|\le\ell(1+|x|)|p-p'|$;
\item\label{compass.3} $(G(t,x,r,p,S)-G(t,x,r',p,S))(r-r')\le\mu|r-r'|^2$;
\item\label{compass.4} for each $R>0$ there is a modulus of continuity $\omega_R$ for which
\[
G(t,x,r,\alpha(x-y),S)-G(t,y,r,\alpha(x-y),S')\le\omega_R\left(\alpha|x-y|^2+|x-y|\right)
\]
for any $\alpha>0$, $t\in[0,T]$, $x,y\in B_R(0)$ and $S,S'\in\Sds^N$ such that
\begin{equation}\label{eq:cilcond}
-3\alpha
\begin{pmatrix}
I_N&0\\0&I_N
\end{pmatrix}
\le
\begin{pmatrix}
S&0\\0&-S'
\end{pmatrix}
\le3\alpha
\begin{pmatrix}
I_N&-I_N\\-I_N&I_N
\end{pmatrix},
\end{equation}
where $I_N$ is the $N\times N$ identity matrix.
\end{myenum}
\end{assum}

For this problem holds a comparison result, which is a parabolic version of \cite[Theorem 7.4]{ishii89}.

\begin{theo}\label{comptheo}
Let $u$ and $v$ be respectively a sub and a supersolution to \eqref{eq:compparpde} under \cref{compass}. If $u(T,x)\le v(T,x)$ for any $x\in\Rds^N$ and both $u$ and $v$ have polynomial growth, then $u\le v$ on $(0,T]\times\Rds^N$.
\end{theo}
\proof
First of all let $\gamma$ be a constant that we will fix later. If we define $u_\gamma(t,x):=e^{\gamma t}u(t,x)$ and $v_\gamma(t,x):=e^{\gamma t}v(t,x)$ we have that they are a sub and a supersolution, respectively, to
\begin{equation}\label{eq:comptheo1}
\partial_tu(t,x)-\gamma u(t,x)+G_\gamma\left(t,x,u,\nabla_xu,D^2_xu\right)=0,
\end{equation}
where $G_\gamma(t,x,r,p,S):=e^{\gamma t}G\left(t,x,e^{-\gamma t}r,e^{-\gamma t}p,e^{-\gamma t}S\right)$. It is a simple verification that $G_\gamma$ still satisfies \cref{compass} and in particular
\begin{equation}\label{eq:comptheo2}
(r-r')(G_\gamma(t,x,r,p,S)-\gamma r-G_\gamma(t,x,r',p,S)+\gamma r')\le(\mu-\gamma)|r-r'|^2.
\end{equation}
Since we assumed that $u$ and $v$ have polynomial growth, there is a constant $k\ge2$ such that, for any $t\in[0,T]$ and $\varepsilon>0$,
\[
\lim_{|x|\to\infty}u_\gamma(t,x)-v_\gamma(t,x)-2\varepsilon(1+|x|)^k=-\infty,
\]
thus we set $\phi(x):=(1+|x|)^k$, $u_{\gamma,\varepsilon}:=u_\gamma-\varepsilon\phi$, $v_{\gamma,\varepsilon}:=v_\gamma+\varepsilon\phi$ and we get
\begin{equation}\label{eq:comptheo3}
\lim_{|x|\to\infty}u_{\gamma,\varepsilon}(t,x)-v_{\gamma,\varepsilon}(t,x)=-\infty.
\end{equation}
Notice that
\begin{align*}
0\le&\partial_tu_\gamma(t,x)-\gamma u_\gamma(t,x)+G_\gamma\left(t,x,u_\gamma,\nabla_xu_\gamma,D^2_xu_\gamma\right)\\
=&\partial_tu_{\gamma,\varepsilon}(t,x)-\gamma(u_{\gamma,\varepsilon}(t,x)+\varepsilon\phi(x))\\
&+G_\gamma\left(t,x,u_{\gamma,\varepsilon}+\varepsilon\phi(x),\nabla_xu_{\gamma,\varepsilon}+\varepsilon\nabla\phi(x),D^2_xu_{\gamma,\varepsilon}+\varepsilon D^2\phi(x)\right)\\
\le&\partial_tu_{\gamma,\varepsilon}(t,x)-\gamma u_{\gamma,\varepsilon}(t,x)-(\gamma-\mu)\varepsilon\phi(x)\\
&+G_\gamma\left(t,x,u_{\gamma,\varepsilon},\nabla_xu_{\gamma,\varepsilon}+\varepsilon\nabla\phi(x),D^2_xu_{\gamma,\varepsilon}+\varepsilon D^2\phi(x)\right)
\end{align*}
and since $|\nabla\phi(x)|$ and $\left|D^2\phi(x)\right|$ are polynomials of degree $k-1$ and $k-2$, respectively, \cref{compass.1,compass.2} in \cref{compass} yield that there exists a constant $C>0$ such that
\begin{align*}
0\le&\partial_tu_{\gamma,\varepsilon}(t,x)-\gamma u_{\gamma,\varepsilon}(t,x)-(\gamma-\mu)\varepsilon\phi(x)\\
&+G_\gamma\left(t,x,u_{\gamma,\varepsilon},\nabla_xu_{\gamma,\varepsilon},D^2_xu_{\gamma,\varepsilon}\right)+C\varepsilon\phi(x).
\end{align*}
Similarly we can prove that
\begin{align*}
0\ge&\partial_tv_{\gamma,\varepsilon}(t,x)-\gamma v_{\gamma,\varepsilon}(t,x)+(\gamma-\mu)\varepsilon\phi(x)\\
&+G_\gamma\left(t,x,v_{\gamma,\varepsilon},\nabla_xv_{\gamma,\varepsilon},D^2_xv_{\gamma,\varepsilon}\right)-C\varepsilon\phi(x).
\end{align*}
Choosing $\gamma=C+\mu>0$ we get that $u_{\gamma,\varepsilon}$ and $v_{\gamma,\varepsilon}$ are respectively a sub and a supersolution to \eqref{eq:comptheo1} and, thanks to \eqref{eq:comptheo2}, $r\mapsto G_\gamma(t,x,r,p,S)-\gamma r$ is decreasing. Finally observe that by \eqref{eq:comptheo3} there is a $R>0$ such that $u_{\gamma,\varepsilon}<v_{\gamma,\varepsilon}$ on $(0,T]\times\Rds^N\setminus B_R(0)$ and that the PDE \eqref{eq:comptheo1} on $(0,T]\times B_R(0)$ satisfies the assumptions of \cite[Theorem 8.2]{userguide}, therefore $u_{\gamma,\varepsilon}\le v_{\gamma,\varepsilon}$ on $(0,T]\times\Rds^N$. The arbitrariness of $\varepsilon$ then implies that $u\le v$ on $(0,T]\times\Rds^N$.
\endproof

It is easy to see that \cref{compass.1,compass.3,compass.2} in \cref{compass} hold true for the operator $F$ in \cref{probass}, therefore, if we can show that it also satisfies \cref{compass.4}, then a comparison result holds for \cref{parprob}. In order to do so let $S,S'\in\Sds^N$ be such that \eqref{eq:cilcond} holds, then, for any $(b,\sigma)\in K_F$, $t\in[0,T]$ and $x\in\Rds^N$, we have
\[
\left(\sigma^\dag S\sigma\right)(t,x)-\left(\sigma^\dag S'\sigma\right)(t,y)\le3\alpha(\sigma(t,x)-\sigma(t,y))^\dag(\sigma(t,x)-\sigma(t,y)),
\]
hence
\begin{align*}
\tr\!\left(\sigma^\dag S\sigma\right)\!(t,x)\!-\!\tr\!\left(\sigma^\dag S'\sigma\right)\!(t,y)\!\le&3\alpha\!\tr\!\left((\sigma(t,x)\!-\!\sigma(t,y))^\dag(\sigma(t,x)\!-\!\sigma(t,y))\right)\\
=&3\alpha|\sigma(t,x)-\sigma(t,y)|^2=3\alpha\ell^2|x-y|^2
\end{align*}
and consequently the identity $\tr\left(AB^\dag\right)=\tr\left(B^\dag A\right)$ yields
\[
\left\langle\sigma\sigma^\dag(t,x),S\right\rangle-\left\langle\sigma\sigma^\dag(t,y),S'\right\rangle\le3\alpha\ell^2|x-y|^2.
\]
From this we get that, for any $\alpha>0$, $t\in[0,T]$, $x,y\in\Rds^N$ and $S,S'\in\Sds^N$ as in \eqref{eq:cilcond}, there exists a constant $C$ such that
\begin{multline*}
F(t,x,r,\alpha(x-y),S)-F(t,y,r,\alpha(x-y),S')\\
\begin{aligned}
\le&\sup_{(b,\sigma)\in K_F}
\begin{aligned}[t]
&\left|\left\langle\sigma\sigma^\dag(t,x),S\right\rangle-\left\langle\sigma\sigma^\dag(t,y),S'\right\rangle\right.\\
&+\alpha(x-y)^\dag b(t,x)-\alpha(x-y)^\dag b(t,y)\\
&\left.+f_{(b,\sigma)}\left(t,x,r,\alpha(x-y)^\dag\sigma\right)-f_{(b,\sigma)}\left(t,y,r,\alpha(x-y)^\dag\sigma\right)\right|
\end{aligned}\\
\le&C\left(\alpha|x-y|^2+|x-y|\right),
\end{aligned}
\end{multline*}
thus $F$ satisfies \cref{compass} and by \cref{comptheo} the next comparison result follows.

\comppar*

\end{appendices}

\emergencystretch=2em
\printbibliography[heading=bibintoc]
\end{document}